\documentclass[12pt]{article}
\usepackage{latexsym,amsfonts,amssymb}
\usepackage[dvips]{color}
\usepackage{colordvi,multicol}
\usepackage{amsmath}

\newtheorem{thm}{Theorem}[section]
\newtheorem{cor}[thm]{Corollary}
\newtheorem{lem}[thm]{Lemma}
\newtheorem{pro}[thm]{Proposition}

\newtheorem{rem}[thm]{Remark}

\newtheorem{con}[thm]{Conjecture}

\date{}
\begin{document}

\title{\bf Some New Gaussian Product
Inequalities}
 \author{Oliver Russell and Wei Sun\\ \\ \\
  {\small Department of Mathematics and Statistics}\\
    {\small Concordia University, Canada}\\ \\
{\small o\_russel@live.concordia.ca,\ \ \ \ wei.sun@concordia.ca}}

\maketitle

\begin{abstract}

\noindent The Gaussian product inequality is a long-standing conjecture. In this paper, we investigate the three-dimensional inequality $E[X_1^{2}X_2^{2m_2}X_3^{2m_3}]\ge E[X_1^{2}]E[X_2^{2m_2}]E[X_3^{2m_3}]$ for any centered Gaussian random vector $(X_1,X_2,X_3)$ and $m_2,m_3\in\mathbb{N}$.
First, we show that this inequality is implied by a combinatorial inequality. The combinatorial inequality can be verified directly for small values of $m_2$ and arbitrary $m_3$. Hence the corresponding cases of the three-dimensional inequality are proved. Second, we show that the three-dimensional inequality is equivalent to an improved Cauchy-Schwarz inequality. This observation leads us to derive some novel moment inequalities for bivariate Gaussian random variables.
\end{abstract}



\section{Introduction}

Correlation is one of the most basic themes of science. Inequalities involving correlations among Gaussian random variables have displayed  importance
and power in both mathematics and physics. In this paper, we study the Gaussian
product inequality (GPI) conjecture. The general form of this conjecture (see Li and Wei \cite{LW12}) says that for any nonnegative real numbers $\alpha_j$, $j=1,\ldots,{n}$, and  any ${n}$-dimensional real-valued centered Gaussian random vector $(X_1,\dots,X_{n})$,
\begin{eqnarray}\label{LW-inequ}
E \left[\prod_{j=1}^{n}|X_j|^{\alpha_j}\right]\geq \prod_{j=1}^{n}E[|X_j|^{\alpha_j}].
\end{eqnarray}
If $\alpha_j=2m$ for $j=1,\ldots,{n}$ and $m\in \mathbb{N}$, then (\ref{LW-inequ}) becomes the Gaussian moment inequality:
\begin{eqnarray}\label{GPC-inequ}
E\left[\prod_{j=1}^{n}X_j^{2m}\right]\geq \prod_{j=1}^{n}E[X_j^{2m}].
\end{eqnarray}
The GPI conjecture is still open. It is known that (cf. Malicet et al. \cite{MNPP16}) if (\ref{GPC-inequ}) holds for any $m\in \mathbb{N}$ then
the `real linear polarization constant' conjecture raised by Ben\'{\i}tem  et al. \cite{BST98} is true.

So far no universal method has been proposed for proving the GPI conjecture; however, several special cases have been solved.
In \cite{Fr08}, Frenkel used algebraic methods to prove (\ref{GPC-inequ}) for the case that $m=1$. In \cite{We14}, Wei used
integral representations to prove a stronger version of (\ref{LW-inequ}) for $\alpha_j\in (-1,0)$ as follows:
\begin{eqnarray*}
E \left[\prod_{j=1}^{n}|X_j|^{\alpha_j}\right]\geq
E \left[\prod_{j=1}^k|X_j|^{\alpha_j}\right]E\left[\prod_{j=k+1}^{n}|X_j|^{\alpha_j}\right],\ \ \ \ \forall 1\le k\le n-1.
\end{eqnarray*}
Malicet et al. \cite{MNPP16} gave a GPI involving Hermite polynomials, which provides a substantial generalization as well as a new analytical proof of Frenkel \cite[Theorem 2.1]{Fr08}.

In \cite[Theorem 3.2]{LHS}, Lan et al. used the Gaussian hypergeometric functions to prove the following inequality: for any $m_1,m_2\in\mathbb{N}$ and any centered Gaussian random vector $(X_1,X_2,X_3)$,
\begin{equation}\label{jan100}
E[X_1^{2m_1}X_2^{2m_2}X_3^{2m_2}]\ge E[X_1^{2m_1}]E[X_2^{2m_2}]E[X_3^{2m_2}].
\end{equation}
Note that the assumption $X_2 $ and $X_3$ have the same exponent $m_2$ plays an essential role in the proof of \cite[Theorem 3.2]{LHS}. A natural question is whether the three-dimensional GPI still holds when the exponents of $X_1,X_2,X_3$ are all different. Surprisingly, this question is very difficult, which has motivated us to consider the following special case of the GPI conjecture.
\begin{con}\label{con2} Let $m_2,m_3\in\mathbb{N}$. For any centered Gaussian random vector $(X_1,X_2,X_3)$,
\begin{equation}\label{jan10}
E[X_1^{2}X_2^{2m_2}X_3^{2m_3}]\ge E[X_1^{2}]E[X_2^{2m_2}]E[X_3^{2m_3}].
\end{equation}
The equality holds if and only if $X_1,X_2,X_3$ are independent.
\end{con}

In Section 3, we will show that Conjecture \ref{con2} is implied by a combinatorial inequality conjecture. The combinatorial inequality can be verified directly for small values of $m_2$ and arbitrary $m_3$. Hence the corresponding cases of Conjecture \ref{con2} are proved (see Theorem \ref{9900} below). In Section 4, we will show that Conjecture \ref{con2} is equivalent to an improved Cauchy-Schwarz inequality. This observation leads us to derive some novel moment inequalities for bivariate Gaussian random variables (see Theorem \ref{thm3133} below). In the next section, we first present some preliminary results, which have independent interest.

\section{Preliminary results}\setcounter{equation}{0}

Throughout this paper, any Gaussian random variable is assumed to be real-valued and non-zero.

\begin{lem}\label{thm2} Let $n\ge 3$ and $m_1,\dots,m_n\in \mathbb{N}$. If for any centered Gaussian random vector $(Y_1,\dots,Y_n)$ with $Y_n=\alpha_1Y_1+\cdots+\alpha_{n-1}Y_{n-1}$ for some constants
$\alpha_1,\dots,\alpha_{n-1}$,
\begin{equation}\label{eqnYn2k}
    E \left[\left\{\prod_{j=1}^{n-1}Y_j^{2m_j}\right\}Y_n^{2k}\right]\ge \left\{\prod_{j=1}^{n-1}E[Y_j^{2m_j}]\right\}E[Y_n^{2k}],\ \ \ \ 0\leq k\leq m_n,
\end{equation}
then for any centered Gaussian random vector $(X_1,\dots,X_n)$,
\begin{equation}\label{14AA}
    E \left[\prod_{j=1}^{n}X_j^{2m_j}\right]\ge \prod_{j=1}^{n}E[X_j^{2m_j}].
\end{equation}
Additionally, if inequality (\ref{eqnYn2k}) is strict when $k=m_n$, then the equality sign of (\ref{14AA}) holds only if $X_n$ is independent of $X_1,\dots,X_{n-1}$.
\end{lem}

\noindent{\bf Proof.}\ \ Let $(X_1,\dots,X_n)$ be a centered Gaussian random vector.
Define
 $$Z_0=E[X_n| X_1,\dots,X_{n-1}], \quad Z_1=X_n-Z_0.$$
Then,
\begin{equation}\label{eqn3Z2k}
    X_n^{2{m_n}}=(Z_0+Z_1)^{2{m_n}}=\sum_{i=0}^{2{m_n}}\binom{2{m_n}}{i} Z_0^{2{m_n}-i} Z_1^i.
\end{equation}
Note that  $Z_1$ is independent of $X_1,\dots,X_{n-1}$. Hence
\begin{equation}\label{eqn3Z0Z1}
  E[Z_0^{2{m_n}-i} Z_1^{i}|  X_1,\dots,X_{n-1}] = Z_0^{2{m_n}-i}E[Z_1^{i}],
\end{equation}
which is equal to zero if $i$ is an odd number.

By (\ref{eqn3Z2k}) and (\ref{eqn3Z0Z1}), we get
\begin{equation}\label{eqn-20EZ|XY}
    E[X_n^{2{m_n}}|  X_1,\dots,X_{n-1}]= \sum_{i=0}^{{m_n}}\binom{2{m_n}}{2i} Z_0^{2{m_n}-2i}E[Z_1^{2i}].
\end{equation}
Note that $Z_0=\alpha_1X_1+\cdots+\alpha_{n-1}X_{n-1}$ for some constants
$\alpha_1,\dots,\alpha_{n-1}$.
Then, it follows from (\ref{eqnYn2k})  that
\begin{equation}\label{eqn3XYZ0}
    E\left[\left(\prod_{j=1}^{n-1}X_j^{2m_j}\right)Z_0^{2{m_n}-2i}\right]
    \geq \left(\prod_{j=1}^{n-1}E[X_j^{2m_j}]\right) E[ Z_0^{2{m_n}-2i}].
\end{equation}
Thus, by (\ref{eqn3Z2k}), (\ref{eqn-20EZ|XY}) and (\ref{eqn3XYZ0}), we obtain  that
\begin{eqnarray}\label{eqn-20E(XYZ1)}
E \left[\prod_{j=1}^{n}X_j^{2m_j}\right]
&=& E \left[\left(\prod_{j=1}^{n-1}X_j^{2m_j}\right)[E[X_n^{2m_n}| X_1,\dots,X_{n-1}]
    \right]\nonumber\\
&=& \sum_{i=0}^{{m_n}}\binom{2{m_n}}{2i}
E\left[\left(\prod_{j=1}^{n-1}X_j^{2m_j}\right)Z_0^{2{m_n}-2i}\right]E[Z_1^{2i}] \nonumber\\
&\geq & \sum_{i=0}^{{m_n}}\binom{2{m_n}}{2i} \left(\prod_{j=1}^{n-1}E[X_j^{2m_j}]\right) E[ Z_0^{2{m_n}-2i}]
E[Z_1^{2i}] \nonumber  \\
&=& \left(\prod_{j=1}^{n-1}E[X_j^{2m_j}]\right)\sum_{i=0}^{{m_n}}\binom{2{m_n}}{2i}
    E[Z_0^{2{m_n}-2i} Z_1^{2i}]   \nonumber\\
&=&\prod_{j=1}^{n}E[X_j^{2m_j}].
\end{eqnarray}

Now, suppose inequality (\ref{eqnYn2k}) is strict when $k=m_n$ and $X_n$ is not independent of $X_1,\dots,X_{n-1}$. Then, $Z_0=E[X_n|X_1,\dots,X_{n-1}]\neq 0$ and
$$
E\left[\left(\prod_{j=1}^{n-1}X_j^{2m_j}\right)Z_0^{2{m_n}}\right] > \left(\prod_{j=1}^{n-1}E[X_j^{2m_j}]\right) E[ Z_0^{2{m_n}}].
$$
Thus, by (\ref{eqn-20E(XYZ1)}), we get
$$
E \left[\prod_{j=1}^{n}X_j^{2m_j}\right] > \prod_{j=1}^{n}E[X_j^{2m_j}].
$$
Therefore, the proof is complete. \hfill\fbox

\begin{cor}\label{corthm2} The following two claims are equivalent.

\noindent Claim I:\ \ For any $n\ge 3$, $m_1,\dots,m_n\in \mathbb{N}$, and centered Gaussian random vector $(X_1,\dots,X_n)$ with $X_n=\alpha_1X_1+\cdots+\alpha_{n-1}X_{n-1}$ for some constants
$\alpha_1,\dots,\alpha_{n-1}$,
$$
E \left[\prod_{j=1}^{n}X_j^{2m_j}\right]> \prod_{j=1}^{n}E[X_j^{2m_j}].
$$

\noindent Claim II:\ \ For any $n\ge 3$, $m_1,\dots,m_n\in \mathbb{N}$, and centered Gaussian random vector $(X_1,\dots,X_n)$,
$$
E \left[\prod_{j=1}^{n}X_j^{2m_j}\right]\ge \prod_{j=1}^{n}E[X_j^{2m_j}],
$$
and the equality holds if and only if $X_1,\dots,X_n$ are independent.
\end{cor}

\noindent{\bf Proof.}\ \ Obviously, Claim II implies Claim I. The assertion that Claim I implies Claim II is a direct consequence of the proof of Lemma \ref{thm2}, where symmetry ensures the equality is equivalent to the independence of $X_1,\dots,X_n$.\hfill\fbox

\begin{lem}\label{lem3}
Let $(X_1,\dots,X_n)$ be a centered Gaussian random vector such that $E[X_iX_j]\ge 0$ for any $i\not=j$. Then,
\begin{eqnarray}\label{24vbnm}
E \left[\prod_{j=1}^{n}X_j^{2m_j}\right]\geq
E \left[\prod_{j=1}^kX_j^{2m_j}\right]E\left[\prod_{j=k+1}^{n}X_j^{2m_j}\right],\ \ \ \ \forall 1\le k\le n-1,
\end{eqnarray}
and
\begin{equation}\label{24V}
E \left[\prod_{j=1}^{n}X_j^{2m_j}\right]\geq\prod_{j=1}^{n}E[X_j^{2m_j}].
\end{equation}
\end{lem}

\noindent{\bf Proof.}\ \ By induction, (\ref{24V}) is a direct consequence of (\ref{24vbnm}). In the following, we prove (\ref{24vbnm}).

We assume without loss of generality that ${\rm Var}(X_i)=1$ for $1\le i\le n$. Denote by $\Lambda$ the covariance matrix of $(X_1,\dots,X_n)$. Let $1\le k\le n-1$. Define
$$
{\tilde\Lambda}_{pq}=\left\{\begin{array}{ll}
\Lambda_{pq},\ \ & {\rm if}\  1\le p,q\le k\ {\rm or}\ k+1\le p,q\le n,\\
0,\ \ & {\rm otherwise}.
\end{array}
\right.
$$
Then, ${\tilde\Lambda}:=({\tilde\Lambda}_{pq})_{1\le p,q\le n}$ is a covariance matrix. Let $(Y_1,\dots,Y_n)$ be a Gaussian random vector with covariance matrix ${\tilde\Lambda}$. Then, $(Y_1,\dots,Y_k)$ and $(X_1,\dots,X_k)$ have the same distribution, $(Y_{k+1},\dots,Y_n)$ and $(X_{k+1},\dots,X_n)$ have the same distribution, and $(Y_1,\dots,Y_k)$ and $(Y_{k+1},\dots,Y_n)$ are independent.

Denote by $T$ the set of all symmetric matrices $s=(s_{kl})_{1\le k,l\le n}$ satisfying $s_{kl}\in \mathbb{N}\cup\{0\}$ and
$$
s_{kk}+\sum_{l=1}^ns_{kl}=2m_k,\ \ \ \ 1\le k\le n.
$$
By \cite[Formula (44)]{SL} (see \cite{Song2017} for the detailed proof), we get
\begin{eqnarray*}
E \left[\prod_{j=1}^{n}X_j^{2m_j}\right]
&=&\left(\prod_{j=1}^n(2m_j)!\right)\sum_{s\in T}\frac{\prod_{k<l}\Lambda_{kl}^{s_{kl}}}{2^{\sum_{k=1}^ns_{kk}}\left(\prod_{k=1}^ns_{kk}!\right)
\left(\prod_{k<l}s_{kl}!\right)}\\
&\ge&\left(\prod_{j=1}^n(2m_j)!\right)\sum_{s\in T}\frac{\prod_{k<l}{{\tilde\Lambda}_{kl}}^{s_{kl}}}{2^{\sum_{k=1}^ns_{kk}}\left(\prod_{k=1}^ns_{kk}!\right)
\left(\prod_{k<l}s_{kl}!\right)}\\
&=&E \left[\prod_{j=1}^{n}Y_j^{2m_j}\right]\\
&=&E \left[\prod_{j=1}^kY_j^{2m_j}\right]E\left[\prod_{j=k+1}^{n}Y_j^{2m_j}\right]\\
&=&E \left[\prod_{j=1}^kX_j^{2m_j}\right]E\left[\prod_{j=k+1}^{n}X_j^{2m_j}\right].
\end{eqnarray*}
The proof is complete.\hfill\fbox

\begin{rem}\label{remn} (i) The GPI conjecture is especially difficult when some components of the Gaussian random vector are negatively correlated.
It is known that if assuming positive correlations then one can verify special cases of the conjecture (cf. e.g., \cite[page 1060]{We14}). Recently, Genest and Ouimet \cite{O} showed that (\ref{24V}) holds if there exists a matrix $C\in[0,\infty)^{n\times n}$ such that $(X_1,\dots,X_n)=(Z_1,\dots,Z_n)C$ in law, where $(Z_1,\dots,Z_n)$
is an $n$-dimensional standard Gaussian random vector. Furthermore, Edelmann et al. \cite{Edel} used a different method to extend (\ref{24vbnm}) to the multivariate gamma distribution.

\noindent (ii) By replacing $X_2$ with $-X_2$ if necessary, Lemma \ref{lem3} implies the well-known two-dimensional GPI:
for any centered bivariate Gaussian random variables $(X_1,X_2)$,
$$
E[X_1^{2m_1}X_2^{2m_2}]\ge E[X_1^{2m_1}]E[X_2^{2m_2}].
$$
\end{rem}

\begin{lem}\label{lem33} Let $m_2,m_3\in \mathbb{N}$ and $U_1$, $U_2$ be independent standard Gaussian random variables. Then, Conjecture \ref{con2} is equivalent to the following claim:

\noindent Claim A:\ \ For any $x_2\in\mathbb{R}$,
\begin{eqnarray}\label{22}
&&(E[U_2(x_2U_1+U_2)^{2m_2}U_1^{2m_3+1}])^2\nonumber\\
&<&\{E[(x_2U_1+U_2)^{2m_2}U_1^{2m_3+2}]-(2m_2-1)!!(2m_3-1)!!(1+x_2^2)^{m_2}\}\nonumber\\
&&\cdot\{E[U_2^{2}(x_2U_1+U_2)^{2m_2}U_1^{2m_3}]-(2m_2-1)!!(2m_3-1)!!(1+x_2^2)^{m_2}\}.
\end{eqnarray}
\end{lem}

\noindent{\bf Proof.}\ \ By Lemma \ref{thm2}, Conjecture \ref{con2} is equivalent to the following inequality:
\begin{eqnarray}\label{1}
&&E[(x_1U_1+U_2)^{2}(x_2U_1+U_2)^{2m_2}U_1^{2m_3}]\nonumber\\
&>&E[(x_1U_1+U_2)^{2}]E[(x_2U_1+U_2)^{2m_2}]E[U_1^{2m_3}],\ \ \ \ \forall x_1,x_2\in\mathbb{R}.\ \ \ \
\end{eqnarray}

We have
\begin{eqnarray*}
&&\ \ (\ref{1})\ {\rm holds}\\
&\Leftrightarrow&\ \ x_1^2E[(x_2U_1+U_2)^{2m_2}U_1^{2m_3+2}]+E[U_2^{2}(x_2U_1+U_2)^{2m_2}U_1^{2m_3}]\\
&&\ \ +2x_1E[U_2(x_2U_1+U_2)^{2m_2}U_1^{2m_3+1}]\nonumber\\
&&>(2m_2-1)!!(2m_3-1)!!(1+x_1^2)(1+x_2^2)^{m_2},\ \ \ \ \forall x_1,x_2\in\mathbb{R}\\
&\Leftrightarrow&\ \ x_1^2\{E[(x_2U_1+U_2)^{2m_2}U_1^{2m_3+2}]-(2m_2-1)!!(2m_3-1)!!(1+x_2^2)^{m_2}\}\\
&&\ \ +2x_1E[U_2(x_2U_1+U_2)^{2m_2}U_1^{2m_3+1}]\nonumber\\
&&\ \ +\{E[U_2^{2}(x_2U_1+U_2)^{2m_2}U_1^{2m_3}]-(2m_2-1)!!(2m_3-1)!!(1+x_2^2)^{m_2}\}\\
&&>0,\ \ \ \ \forall x_1,x_2\in\mathbb{R}.
\end{eqnarray*}
Note that
\begin{eqnarray*}
&&E[(x_2U_1+U_2)^{2m_2}U_1^{2m_3+2}]-(2m_2-1)!!(2m_3-1)!!(1+x_2^2)^{m_2}\\
&\ge&E[(x_2U_1+U_2)^{2m_2}]E[U_1^{2m_3+2}]-(2m_2-1)!!(2m_3-1)!!(1+x_2^2)^{m_2}\\
&=&(2m_2-1)!!(2m_3+1)!!(1+x_2^2)^{m_2}-(2m_2-1)!!(2m_3-1)!!(1+x_2^2)^{m_2}\\
&=&2m_3(2m_2-1)!!(2m_3-1)!!(1+x_2^2)^{m_2}\\
&>&0,
\end{eqnarray*}
and  Lemma \ref{lem3} implies that
\begin{eqnarray*}
&&E[U_2^{2}(x_2U_1+U_2)^{2m_2}U_1^{2m_3}]-(2m_2-1)!!(2m_3-1)!!(1+x_2^2)^{m_2}\\
&=&E[U_2^{2}(|x_2|U_1+U_2)^{2m_2}U_1^{2m_3}]-(2m_2-1)!!(2m_3-1)!!(1+x_2^2)^{m_2}\\
&\ge&(2m_2-1)!!(2m_3-1)!!(1+x_2^2)^{m_2}-(2m_2-1)!!(2m_3-1)!!(1+x_2^2)^{m_2}\\
&=&0.
\end{eqnarray*}
Hence, (\ref{1}) holds if and only if (\ref{22}) holds for any $x_2\in\mathbb{R}$.
The proof is complete.\hfill\fbox

\section{Combinatorial inequality and proof of GPI for small $m_2$}\setcounter{equation}{0}

In this section, we propose a combinatorial inequality conjecture (see Conjecture \ref{con1} below) and show that it implies Conjecture 1.1. Although we have not been able to completely resolve this combinatorial inequality conjecture, it can be verified directly for small values of $m_2$ and arbitrary
$m_3$. Thereby, the corresponding cases of Conjecture 1.1 are proved (see Theorem \ref{9900} below).

Let $m_2,m_3\in \mathbb{N}$. For $0\le p\le 2m_2$, define
\begin{eqnarray*}
c^{m_2,m_3}_p&=&\sum_{i=\max(0,p-m_2)}^{\min(m_2,p)}\left[\frac{2^{2i}\left(m_3+i+\frac{1}{2}\right)\cdots\frac{1}{2}}{(2i)!(m_2-i)!}
-\frac{\left(m_3-\frac{1}{2}\right)\cdots\frac{1}{2}}{2(i!)(m_2-i)!}\right]\nonumber\\
&&\ \ \ \ \cdot\left[\frac{2^{2p-2i-1}(2m_2+2i-2p+1)\left(m_3+p-i-\frac{1}{2}\right)\cdots\frac{1}{2}}{(2p-2i)!(m_2+i-p)!}-\frac{\left(m_3-\frac{1}{2}\right)
\cdots\frac{1}{2}}{2[(p-i)!](m_2+i-p)!}\right]\nonumber\\
&&-\sum_{i=\max(0,p-m_2)}^{\min(m_2,p)-1}\frac{2^{2p}\left(m_3+i+\frac{1}{2}\right)\cdots\frac{1}{2}
\left(m_3+p-i-\frac{1}{2}\right)\cdots\frac{1}{2}}
{(2i+1)!(2p-2i-1)!(m_2-i-1)!(m_2+i-p)!}.
\end{eqnarray*}

\begin{lem}\label{7878}  Let $m_2,m_3\in \mathbb{N}$. Then, we have

\noindent (i) $c^{m_2,m_3}_0,c^{m_2,m_3}_2>0$; $c^{m_2,m_3}_1\geq0$ if $m_3\le m_2-1$.

\noindent (ii) $(c^{m_2,m_3}_1)^2< 4c^{m_2,m_3}_0c^{m_2,m_3}_2$.

\noindent (iii) $c^{m_2,m_3}_{2m_2}>0$.

\end{lem}

\noindent{\bf Proof.}\ \ To simplify notation, we denote $c^{m_2,m_3}_p$ by $c_p$ for $0\le p\le 2m_2$. Define $r:=\left(m_3-\frac{1}{2}\right)\cdots\frac{1}{2}$.

We have
\begin{eqnarray*}
c_{0}&=&\frac{r^2 m_3}{m_2!(m_2-1)!},
\end{eqnarray*}
\begin{eqnarray*}
c_{1}&=&r^2\cdot\frac{m_3}{m_2!}\cdot\frac{2m_2 m_3+m_2-m_3-1}{(m_2-1)!}+r^2\cdot\frac{2m_3^2+4m_3+1}{(m_2-1)!}\cdot\frac{m_2}{m_2!}-r^2\cdot\frac{4\left(m_3+\frac{1}{2}\right)^2 m_2}{(m_2-1)!m_2!}\\
&=&-\frac{r^2 m_3}{m_2!(m_2-1)!}\cdot(m_3+1-m_2),
\end{eqnarray*}
\begin{eqnarray*}
c_{2}&=&r^2\cdot\frac{m_3}{m_2!}\cdot\frac{8(2m_2-3)(m_3+\frac{3}{2})(m_3+\frac{1}{2})-6}{4!(m_2-1)!}\cdot(m_2-1)\\
&&+r^2\cdot\frac{12(2m_3^2+4m_3+1)}{(m_2-1)!}\cdot\frac{2(2m_2-1)\left(m_3+\frac{1}{2}\right)-1}{4!m_2!}\cdot m_2\\
&&+r^2\cdot\frac{16\left(m_3+\frac{5}{2}\right)\left(m_3+\frac{3}{2}\right)\left(m_3+\frac{1}{2}\right)-6}{4!(m_2-1)!}\cdot\frac{m_2(m_2-1)}{m_2!}\\
&&-r^2\cdot\frac{128\left(m_3+\frac{1}{2}\right)^2\left(m_3+\frac{3}{2}\right)}{4!m_2!(m_2-1)!}\cdot m_2(m_2-1)\\
&=&\frac{r^2 m_3}{m_2!(m_2-1)!}\cdot(m_2m_3^2+m_3^2+m_2^2m_3+m_2m_3+2m_3+m_2^2+1-m_2)\\
&>&0,
\end{eqnarray*}
\begin{eqnarray*}
c_{2m_2}&=&\left[\frac{2^{2m_2}\left(m_3+m_2+\frac{1}{2}\right)\cdots\frac{1}{2}}{(2m_2)!}-\frac{\left(m_3-\frac{1}{2}\right)\cdots\frac{1}{2}}{2(m_2!)}\right]\\
&&\cdot\left[\frac{2^{2m_2-1}\left(m_3+m_2-\frac{1}{2}\right)\cdots\frac{1}{2}}{(2m_2)!}-\frac{\left(m_3-\frac{1}{2}\right)\cdots\frac{1}{2}}{2(m_2!)}
\right]\\
&=&\left[\frac{\left(m_3+m_2+\frac{1}{2}\right)\cdots(m_2+\frac{1}{2})}{m_2!}-\frac{\left(m_3-\frac{1}{2}\right)\cdots\frac{1}{2}}{2(m_2!)}\right]\\
&&\cdot\left[\frac{\left(m_3+m_2-\frac{1}{2}\right)\cdots(m_2+\frac{1}{2})}{2(m_2!)}-\frac{\left(m_3-\frac{1}{2}\right)\cdots\frac{1}{2}}{2(m_2!)}
\right]\\
&>&0,
\end{eqnarray*}
and
\begin{eqnarray*}
c_1^2<4c_0c_2&\Leftrightarrow&(m_3+1-m_2)^2<4(m_2m_3^2+m_3^2+m_2^2m_3+m_2m_3+2m_3+m_2^2+1-m_2)\\
&\Leftrightarrow&4m_2m_3^2+3m_3^2+4m_2^2m_3+6m_2m_3+6m_3+3m_2^2+3-2m_2>0,
\end{eqnarray*}
which is clearly true.\hfill\fbox

Now we state the combinatorial inequality conjecture.
\begin{con}\label{con1}
Let $m_2\ge 2$ and $m_3\in \mathbb{N}$. Then, we have
\begin{equation}\label{rtrtrt}
c^{m_2,m_3}_p>0\ {\rm for}\ 3\le p\le 2m_2-1.
\end{equation}
\end{con}

\begin{pro}\label{hjhj}
If Conjecture \ref{con1} is true, then Conjecture 1.1 is also true.
\end{pro}

\noindent{\bf Proof.}\ \ Assume that Conjecture \ref{con1} is true. We we will apply Lemma \ref{lem33} to show that Conjecture 1.1 is also true.

Let $U_1$, $U_2$ be independent standard Gaussian random variables. For $x_2\in\mathbb{R}$, we have
\begin{eqnarray*}
&&(E[U_2(x_2U_1+U_2)^{2m_2}U_1^{2m_3+1}])^2\nonumber\\
&=&\left(\sum_{k\ {\rm odd}}{2m_2\choose k}x_2^k(2m_3+k)!!(2m_2-k)!!\right)^2\nonumber\\
&=&\sum_{k\ {\rm odd}}\sum_{l\ {\rm odd}}x_2^{k+l}{2m_2\choose k}{2m_2\choose l}(2m_3+k)!!(2m_2-k)!!(2m_3+l)!!(2m_2-l)!!,
\end{eqnarray*}
\begin{eqnarray*}
&&E[(x_2U_1+U_2)^{2m_2}U_1^{2m_3+2}]-(2m_2-1)!!(2m_3-1)!!(1+x_2^2)^{m_2}\nonumber\\
&=&\sum_{k\ {\rm even}}x_2^k\left[{2m_2\choose k}(2m_3+k+1)!!(2m_2-k-1)!!-{m_2\choose \frac{k}{2}}(2m_2-1)!!(2m_3-1)!!\right],
\end{eqnarray*}
and
\begin{eqnarray*}
&&E[U_2^{2}(x_2U_1+U_2)^{2m_2}U_1^{2m_3}]-(2m_2-1)!!(2m_3-1)!!(1+x_2^2)^{m_2}\nonumber\\
&=&\sum_{l\ {\rm even}}x_2^l\left[{2m_2\choose l}(2m_3+l-1)!!(2m_2-l+1)!!-{m_2\choose \frac{l}{2}}(2m_2-1)!!(2m_3-1)!!\right].
\end{eqnarray*}
Hence, inequality (\ref{22}) becomes
\begin{eqnarray}\label{12AS}
&&\sum_{k\ {\rm odd}}\sum_{l\ {\rm odd}}x_2^{k+l}{2m_2\choose k}{2m_2\choose l}(2m_3+k)!!(2m_2-k)!!(2m_3+l)!!(2m_2-l)!!\nonumber\\
&<&\sum_{k\ {\rm even}}\sum_{l\ {\rm even}}x_2^{k+l}\left[{2m_2\choose k}(2m_3+k+1)!!(2m_2-k-1)!!-{m_2\choose \frac{k}{2}}(2m_2-1)!!(2m_3-1)!!\right]\nonumber\\
&&\ \ \ \ \cdot\left[{2m_2\choose l}(2m_3+l-1)!!(2m_2-l+1)!!-{m_2\choose \frac{l}{2}}(2m_2-1)!!(2m_3-1)!!\right].
\end{eqnarray}

We have
\begin{eqnarray*}
&&{2m_2\choose k}{2m_2\choose l}(2m_3+k)!!(2m_2-k)!!(2m_3+l)!!(2m_2-l)!!\\
&=&\frac{(2m_2)!}{k!(2m_2-k)!}\cdot\frac{(2m_2)!}{l!(2m_2-l)!}\cdot(2m_2-k)!!(2m_2-l)!!(2m_3+k)!!(2m_3+l)!!\\
&=&\frac{(2m_2)!}{k!2^{m_2-\frac{k+1}{2}}(m_2-\frac{k+1}{2})!}\cdot\frac{(2m_2)!}{l!2^{m_2-\frac{l+1}{2}}(m_2-\frac{l+1}{2})!}\cdot(2m_3+k)!!(2m_3+l)!!\\
&=&\frac{(2m_2)!(2m_2)!}{2^{2m_2-\frac{k+l}{2}-1}k!l!(m_2-\frac{k+1}{2})!(m_2-\frac{l+1}{2})!}\\
&&\cdot2^{m_3+\frac{k+1}{2}}\cdot\left(m_3+\frac{k}{2}\right)\cdots\frac{1}{2}
\cdot2^{m_3+\frac{l+1}{2}}\cdot\left(m_3+\frac{l}{2}\right)\cdots\frac{1}{2}\\
&=&\frac{2^{2m_3-2m_2+k+l+2}(2m_2)!(2m_2)!\left(m_3+\frac{k}{2}\right)\cdots\frac{1}{2}\left(m_3+\frac{l}{2}\right)\cdots\frac{1}{2}}
{k!l!(m_2-\frac{k+1}{2})!(m_2-\frac{l+1}{2})!},
\end{eqnarray*}
\begin{eqnarray*}
&&{2m_2\choose k}(2m_3+k+1)!!(2m_2-k-1)!!-{m_2\choose \frac{k}{2}}(2m_2-1)!!(2m_3-1)!!\\
&=&\frac{(2m_2)!}{k!(2m_2-k)!}\cdot(2m_2-k-1)!!\cdot2^{m_3+\frac{k+2}{2}}\cdot\left(m_3+\frac{k+1}{2}\right)\cdots\frac{1}{2}\\
&&-\frac{m_2!}{(\frac{k}{2})!(m_2-\frac{k}{2})!}\cdot2^{m_2}\cdot\left(m_2-\frac{1}{2}\right)\cdots\frac{1}{2}
\cdot2^{m_3}\cdot\left(m_3-\frac{1}{2}\right)\cdots\frac{1}{2}\\
&=&\frac{(2m_2)!}{k!(m_2-\frac{k}{2})!}\cdot2^{m_3-m_2+k+1}\cdot\left(m_3+\frac{k+1}{2}\right)\cdots\frac{1}{2}\\
&&-\frac{m_2!}{(\frac{k}{2})!(m_2-\frac{k}{2})!}\cdot2^{m_2+m_3}\cdot
\left(m_2-\frac{1}{2}\right)\cdots\frac{1}{2}\left(m_3-\frac{1}{2}\right)\cdots\frac{1}{2},
\end{eqnarray*}
and
\begin{eqnarray*}
&&{2m_2\choose l}(2m_3+l-1)!!(2m_2-l+1)!!-{m_2\choose \frac{l}{2}}(2m_2-1)!!(2m_3-1)!!\\
&=&\frac{(2m_2)!}{l!(2m_2-l)!}\cdot(2m_2-l+1)!!\cdot2^{m_3+\frac{l}{2}}\cdot\left(m_3+\frac{l-1}{2}\right)\cdots\frac{1}{2}\\
&&-\frac{m_2!}{(\frac{l}{2})!(m_2-\frac{l}{2})!}\cdot2^{m_2}\cdot\left(m_2-\frac{1}{2}\right)\cdots\frac{1}{2}
\cdot2^{m_3}\cdot\left(m_3-\frac{1}{2}\right)\cdots\frac{1}{2}\\
&=&\frac{(2m_2)!(2m_2-l+1)}{l!(m_2-\frac{l}{2})!}\cdot2^{m_3-m_2+l}\cdot\left(m_3+\frac{l-1}{2}\right)\cdots\frac{1}{2}\\
&&-\frac{m_2!}{(\frac{l}{2})!(m_2-\frac{l}{2})!}\cdot2^{m_2+m_3}\cdot\left(m_2-\frac{1}{2}\right)\cdots\frac{1}{2}\left(m_3-\frac{1}{2}\right)
\cdots\frac{1}{2}.
\end{eqnarray*}
Then,
\begin{eqnarray*}
&&\left[{2m_2\choose k}(2m_3+k+1)!!(2m_2-k-1)!!-{m_2\choose \frac{k}{2}}(2m_2-1)!!(2m_3-1)!!\right]\\
&&\cdot\left[{2m_2\choose l}(2m_3+l-1)!!(2m_2-l+1)!!-{m_2\choose \frac{l}{2}}(2m_2-1)!!(2m_3-1)!!\right]\\
&=&\left[\frac{(2m_2)!}{k!(m_2-\frac{k}{2})!}\cdot2^{m_3-m_2+k+1}\cdot\left(m_3+\frac{k+1}{2}\right)\cdots\frac{1}{2}\right.\\
&&\left.-\frac{m_2!}{(\frac{k}{2})!(m_2-\frac{k}{2})!}\cdot2^{m_2+m_3}\cdot
\left(m_2-\frac{1}{2}\right)\cdots\frac{1}{2}\left(m_3-\frac{1}{2}\right)\cdots\frac{1}{2}\right]\\
&&\cdot\left[\frac{(2m_2)!(2m_2-l+1)}{l!(m_2-\frac{l}{2})!}\cdot2^{m_3-m_2+l}\cdot\left(m_3+\frac{l-1}{2}\right)\cdots\frac{1}{2}\right.\\
&&\left.-\frac{m_2!}{(\frac{l}{2})!(m_2-\frac{l}{2})!}\cdot2^{m_2+m_3}\cdot\left(m_2-\frac{1}{2}\right)\cdots\frac{1}{2}\left(m_3-\frac{1}{2}\right)
\cdots\frac{1}{2}\right].
\end{eqnarray*}
Hence, (\ref{12AS}) becomes
\begin{eqnarray}\label{popo}
&&\sum_{k\ {\rm odd}}\sum_{l\ {\rm odd}}x_2^{k+l}\cdot\frac{2^{2m_3-2m_2+k+l+2}(2m_2)!(2m_2)!\left(m_3+\frac{k}{2}\right)\cdots\frac{1}{2}\left(m_3+\frac{l}{2}\right)\cdots\frac{1}{2}}
{k!l!(m_2-\frac{k+1}{2})!(m_2-\frac{l+1}{2})!}\nonumber\\
&<&\sum_{k\ {\rm even}}\sum_{l\ {\rm even}}x_2^{k+l}\left[\frac{(2m_2)!}{k!(m_2-\frac{k}{2})!}\cdot2^{m_3-m_2+k+1}\cdot\left(m_3+\frac{k+1}{2}\right)\cdots\frac{1}{2}\right.\nonumber\\
&&\left.-\frac{m_2!}{(\frac{k}{2})!(m_2-\frac{k}{2})!}\cdot2^{m_2+m_3}\cdot
\left(m_2-\frac{1}{2}\right)\cdots\frac{1}{2}\left(m_3-\frac{1}{2}\right)\cdots\frac{1}{2}\right]\nonumber\\
&&\cdot\left[\frac{(2m_2)!(2m_2-l+1)}{l!(m_2-\frac{l}{2})!}\cdot2^{m_3-m_2+l}\cdot\left(m_3+\frac{l-1}{2}\right)\cdots\frac{1}{2}\right.\nonumber\\
&&\left.-\frac{m_2!}{(\frac{l}{2})!(m_2-\frac{l}{2})!}\cdot2^{m_2+m_3}\cdot\left(m_2-\frac{1}{2}\right)\cdots\frac{1}{2}\left(m_3-\frac{1}{2}\right)
\cdots\frac{1}{2}\right].
\end{eqnarray}
Diving both sides of inequality (\ref{popo}) by $2^{m_3-m_2+1}(2m_2)!(2m_2)!$, we get
\begin{eqnarray}\label{15Z}
&&\sum_{k\ {\rm odd}}\sum_{l\ {\rm odd}}x_2^{k+l}\cdot\frac{2^{k+l}\left(m_3+\frac{k}{2}\right)\cdots\frac{1}{2}\left(m_3+\frac{l}{2}\right)\cdots\frac{1}{2}}
{k!l!(m_2-\frac{k+1}{2})!(m_2-\frac{l+1}{2})!}\nonumber\\
&<&\sum_{k\ {\rm even}}\sum_{l\ {\rm even}}x_2^{k+l}\left[\frac{2^{k}\left(m_3+\frac{k+1}{2}\right)\cdots\frac{1}{2}}{k!(m_2-\frac{k}{2})!}
-\frac{\left(m_3-\frac{1}{2}\right)\cdots\frac{1}{2}}{2(\frac{k}{2})!(m_2-\frac{k}{2})!}
\right]\nonumber\\
&&\cdot\left[\frac{2^{l-1}(2m_2-l+1)\left(m_3+\frac{l-1}{2}\right)\cdots\frac{1}{2}}{l!(m_2-\frac{l}{2})!}-\frac{\left(m_3-\frac{1}{2}\right)
\cdots\frac{1}{2}}{2(\frac{l}{2})!(m_2-\frac{l}{2})!}\right].
\end{eqnarray}

Define
$$
x=x_2^2.
$$
We have
\begin{eqnarray*}
&&\sum_{k\ {\rm odd}}\sum_{l\ {\rm odd}}x_2^{k+l}\cdot\frac{2^{k+l}\left(m_3+\frac{k}{2}\right)\cdots\frac{1}{2}\left(m_3+\frac{l}{2}\right)\cdots\frac{1}{2}}
{k!l!(m_2-\frac{k+1}{2})!(m_2-\frac{l+1}{2})!}\nonumber\\
&=&\sum_{i=1}^{m_2}\sum_{j=1}^{m_2}x_2^{2(i+j-1)}\cdot\frac{2^{2(i+j-1)}\left(m_3+i-\frac{1}{2}\right)\cdots\frac{1}{2}
\left(m_3+j-\frac{1}{2}\right)\cdots\frac{1}{2}}
{(2i-1)!(2j-1)!(m_2-i)!(m_2-j)!}\nonumber\\
&=&\sum_{p=1}^{2m_2-1}x^{p}\sum_{i=\max(1,p+1-m_2)}^{\min(m_2,p)}\frac{2^{2p}\left(m_3+i-\frac{1}{2}\right)\cdots\frac{1}{2}
\left(m_3+p-i+\frac{1}{2}\right)\cdots\frac{1}{2}}
{(2i-1)!(2p-2i+1)!(m_2-i)!(m_2+i-p-1)!}\\
&=&\sum_{p=1}^{2m_2-1}x^{p}\sum_{i=\max(0,p-m_2)}^{\min(m_2,p)-1}\frac{2^{2p}\left(m_3+i+\frac{1}{2}\right)\cdots\frac{1}{2}
\left(m_3+p-i-\frac{1}{2}\right)\cdots\frac{1}{2}}
{(2i+1)!(2p-2i-1)!(m_2-i-1)!(m_2+i-p)!},
\end{eqnarray*}
and
\begin{eqnarray*}
&&\sum_{k\ {\rm even}}\sum_{l\ {\rm even}}x_2^{k+l}\left[\frac{2^{k}\left(m_3+\frac{k+1}{2}\right)\cdots\frac{1}{2}}{k!(m_2-\frac{k}{2})!}
-\frac{\left(m_3-\frac{1}{2}\right)\cdots\frac{1}{2}}{2(\frac{k}{2})!(m_2-\frac{k}{2})!}
\right]\nonumber\\
&&\cdot\left[\frac{2^{l-1}(2m_2-l+1)\left(m_3+\frac{l-1}{2}\right)\cdots\frac{1}{2}}{l!(m_2-\frac{l}{2})!}-\frac{\left(m_3-\frac{1}{2}\right)
\cdots\frac{1}{2}}{2(\frac{l}{2})!(m_2-\frac{l}{2})!}\right]\nonumber\\
&=&\sum_{i=0}^{m_2}\sum_{j=0}^{m_2}x_2^{2(i+j)}\left[\frac{2^{2i}\left(m_3+i+\frac{1}{2}\right)\cdots\frac{1}{2}}{(2i)!(m_2-i)!}
-\frac{\left(m_3-\frac{1}{2}\right)\cdots\frac{1}{2}}{2(i!)(m_2-i)!}\right]\\
&&\cdot\left[\frac{2^{2j-1}(2m_2-2j+1)\left(m_3+j-\frac{1}{2}\right)\cdots\frac{1}{2}}{(2j)!(m_2-j)!}-\frac{\left(m_3-\frac{1}{2}\right)
\cdots\frac{1}{2}}{2(j!)(m_2-j)!}\right]\nonumber\\
&=&\sum_{p=0}^{2m_2}x^{p}\sum_{i=\max(0,p-m_2)}^{\min(m_2,p)}\left[\frac{2^{2i}\left(m_3+i+\frac{1}{2}\right)\cdots\frac{1}{2}}{(2i)!(m_2-i)!}
-\frac{\left(m_3-\frac{1}{2}\right)\cdots\frac{1}{2}}{2(i!)(m_2-i)!}\right]\\
&&\cdot\left[\frac{2^{2p-2i-1}(2m_2+2i-2p+1)\left(m_3+p-i-\frac{1}{2}\right)\cdots\frac{1}{2}}{(2p-2i)!(m_2+i-p)!}-\frac{\left(m_3-\frac{1}{2}\right)
\cdots\frac{1}{2}}{2[(p-i)!](m_2+i-p)!}\right].
\end{eqnarray*}
Define $F(x)$ to be the polynomial given by (RHS - LHS) of inequality (\ref{15Z}). Then, we have that
$$
F(x)=\sum_{p=0}^{2m_2}c^{m_2,m_3}_px^{p}.
$$

By  (\ref{rtrtrt}) and Lemma \ref{7878}, we get
$$
F(x)\ge c^{m_2,m_3}_0+c^{m_2,m_3}_1x+c^{m_2,m_3}_2x^2>0,\ \ \ \ \forall x\ge0.
$$
Then, inequality (\ref{15Z}) and hence inequality (\ref{22}) hold. Therefore, the proof is complete by Lemma \ref{lem33}.\hfill\fbox

\begin{lem}\label{78}  Let $m_2\in\{2,3\}$ and $m_3\in \mathbb{N}$. Then, we have
$$
c^{m_2,m_3}_p>0\ {\rm for}\ 3\le p\le 2m_2-1.
$$
\end{lem}

\noindent{\bf Proof.}\ \ To simplify notation, we denote $c^{m_2,m_3}_p$ by $c_p$ for $0\le p\le 2m_2$. Define $r:=\left(m_3-\frac{1}{2}\right)\cdots\frac{1}{2}$.

First, we consider the case $m_2=2$. We only need to show that $c_3>0$. We have
\begin{eqnarray*}
c_{3}&=&r^2\cdot\frac{1}{3}m_3\left(m_3+2\right)\left(2m_3^2+4m_3+1\right)\\
&&+r^2\cdot\frac{1}{6}\left(3m_3+1\right)\left(4m_3^3+18m_3^2+23m_3+6\right)\\
&&-r^2\cdot\frac{16}{9}\left(m_3+\frac{3}{2}\right)^2\left(m_3+\frac{1}{2}\right)^2\\
&=&\frac{r^2 m_3}{18}\cdot\left(16m_3^3+94m_3^2+139m_3+39\right)>0,
\end{eqnarray*}
which completes the proof for this case.

Next, consider the case $m_2=3$. We need to show that
$$
c_3>0,\ c_4>0\ \ {\rm and}\ \ c_5>0.
$$
We have
\begin{eqnarray*}
c_{3}&=&r^2\cdot\frac{1}{540}m_3^2\left(4m_3^2+18m_3+23\right)\\
&&+r^2\cdot\left(m_3^2+2m_3+\frac{1}{2}\right)^2\\
&&+r^2\cdot\frac{1}{12}\left(4m_3^3+18m_3^2+23m_3+6\right)\left(5m_3+2\right)\\
&&+r^2\cdot\frac{1}{180}\left(8m_3^4+64m_3^3+172m_3^2+176m_3+45\right)\\
&&-r^2\cdot\frac{8}{15}\left(m_3+\frac{5}{2}\right)\left(m_3+\frac{3}{2}\right)\left(m_3+\frac{1}{2}\right)^2\\
&&-r^2\cdot\frac{16}{9}\left(m_3+\frac{3}{2}\right)^2\left(m_3+\frac{1}{2}\right)^2\\
&=&\frac{r^2 m_3}{54}\cdot\left(22m_3^3+150m_3^2+245m_3+78\right)>0,
\end{eqnarray*}
\begin{eqnarray*}
c_{4}&=&r^2\cdot\frac{1}{90}m_3\left(m_3^2+2m_3+\frac{1}{2}\right)\left(4m_3^2+18m_3+23\right)\\
&&+r^2\cdot\frac{1}{6}\left(4m_3^3+18m_3^2+23m_3+6\right)\left(m_3^2+2m_3+\frac{1}{2}\right)\\
&&+r^2\cdot\frac{1}{180}\left(8m_3^4+64m_3^3+172m_3^2+176m_3+45\right)\left(5m_3+2\right)\\
&&-r^2\cdot\frac{32}{45}\left(m_3+\frac{5}{2}\right)\left(m_3+\frac{3}{2}\right)^2\left(m_3+\frac{1}{2}\right)^2\\
&=&\frac{r^2 m_3}{180}\cdot\left(40m_3^4+336m_3^3+956m_3^2+1020m_3+273\right)>0,
\end{eqnarray*}
and
\begin{eqnarray*}
c_{5}&=&r^2\cdot\frac{1}{540}m_3\left(4m_3^3+18m_3^2+23m_3+6\right)\left(4m_3^2+18m_3+23\right)\\
&&+r^2\cdot\frac{1}{90}\left(8m_3^4+64m_3^3+172m_3^2+176m_3+45\right)\left(m_3^2+2m_3+\frac{1}{2}\right)\\
&&-r^2\cdot\frac{16}{225}\left(m_3+\frac{5}{2}\right)^2\left(m_3+\frac{3}{2}\right)^2\left(m_3+\frac{1}{2}\right)^2\\
&=&\frac{r^2 m_3}{2700}\cdot\left(128m_3^5+1392m_3^4+5564m_3^3+10164m_3^2+8087m_3+1890\right)>0,
\end{eqnarray*}
which completes the proof for this case.\hfill\fbox

\begin{thm}\label{9900} Let $m_3\in\mathbb{N}$. For any centered Gaussian random vector $(X_1,X_2,X_3)$,
$$
E[X_1^{2}X_2^{4}X_3^{2m_3}]\ge E[X_1^{2}]E[X_2^{4}]E[X_3^{2m_3}],
$$
and
$$
E[X_1^{2}X_2^{6}X_3^{2m_3}]\ge E[X_1^{2}]E[X_2^{6}]E[X_3^{2m_3}].
$$
The equalities hold if and only if $X_1,X_2,X_3$ are independent.
\end{thm}

\noindent{\bf Proof.}\ \ This is a direct consequence of Proposition \ref{hjhj} and Lemma \ref{78}.\hfill\fbox

\section{Improved Cauchy-Schwarz inequality for bivariate Gaussian random variables}\setcounter{equation}{0}

In this section, we continue studying Conjecture \ref{con2}. First, we show that Conjecture \ref{con2} is equivalent to an improved Cauchy-Schwarz inequality.

\begin{lem}\label{lem44} Let $m_2,m_3\in \mathbb{N}$. Then, Conjecture \ref{con2} is equivalent to the following claim:

\noindent Claim B:\ \ Let $(U,V)$ be bivariate Gaussian random variables with $U\sim N(0,1)$ and $V\sim N(0,1)$. Then,
\begin{eqnarray}\label{14DD}
&&(E[ U^{2m_2+1}V^{2m_3+1}])^2+(2m_2-1)!!(2m_3-1)!!\{E[ U^{2m_2}V^{2m_3+2}]+E[V^{2m_2+2}U^{2m_3}]\}\nonumber\\
&\le&E[U^{2m_2}V^{2m_3+2}] E[V^{2m_2+2}U^{2m_3}]\nonumber\\
&&+(2m_2-1)!!(2m_3-1)!!\{2E[ U^{2m_2+1}V^{2m_3+1}] E[UV]+(2m_2-1)!!(2m_3-1)!!(1-(E[UV])^2)\}.\nonumber\\
&&
\end{eqnarray}
The equality holds if and only if $|E[UV]|=1$.
\end{lem}

\noindent{\bf Proof.}\ \  By Lemma \ref{lem33}, we only need to show that inequality (\ref{22}) is equivalent to inequality (\ref{14DD}).

Obviously, (\ref{22}) holds if $x=0$. To prove (\ref{22}), we assume without loss of generality that $x_2> 0$. Define
$$
\theta=\arccos\left(\frac{x_2}{\sqrt{1+x_2^2}}\right).
$$
Then, $\theta\in(0,\pi/2)$ and (\ref{22}) is equivalent to
\begin{eqnarray}\label{Q1}
&&(E[U_2( U_1\cos\theta+U_2\sin\theta)^{2m_2}U_1^{2m_3+1}])^2\nonumber\\
&<&\{E[(U_1\cos\theta+U_2\sin\theta)^{2m_2}U_1^{2m_3+2}]-(2m_2-1)!!(2m_3-1)!!\}\nonumber\\
&&\cdot\{E[U_2^{2}(U_1\cos\theta+U_2\sin\theta)^{2m_2}U_1^{2m_3}]-(2m_2-1)!!(2m_3-1)!!\},\ \ \ \ \forall \theta\in(0,\pi/2).\ \ \ \
\end{eqnarray}
Define
$$
V= U_1\cos\theta+U_2\sin\theta.
$$
Then, we have that
\begin{eqnarray*}
&&\ \ (\ref{Q1})\ {\rm holds}\nonumber\\
&\Leftrightarrow&\ \ (E[U_2V^{2m_2}U_1^{2m_3+1}])^2\nonumber\\
&&<\{E[V^{2m_2}U_1^{2m_3+2}]-(2m_2-1)!!(2m_3-1)!!\}\nonumber\\
&&\ \ \ \ \cdot\{E[U_2^{2}V^{2m_2}U_1^{2m_3}]-(2m_2-1)!!(2m_3-1)!!\}\\
&\Leftrightarrow&\ \ (E[(U_2\sin\theta) V^{2m_2}U_1^{2m_3+1}])^2\nonumber\\
&&<\{E[V^{2m_2}U_1^{2m_3+2}]-(2m_2-1)!!(2m_3-1)!!\}\nonumber\\
&&\ \ \ \ \cdot\{E[(U_2\sin\theta)^{2}V^{2m_2}U_1^{2m_3}]-(2m_2-1)!!(2m_3-1)!!\sin^2\theta\}\\
&\Leftrightarrow&\ \ (E[(V-U_1\cos\theta) V^{2m_2}U_1^{2m_3+1}])^2\nonumber\\
&&<\{E[V^{2m_2}U_1^{2m_3+2}]-(2m_2-1)!!(2m_3-1)!!\}\nonumber\\
&&\ \ \ \ \cdot\{E[(V- U_1\cos\theta)^{2}V^{2m_2}U_1^{2m_3}]-(2m_2-1)!!(2m_3-1)!!\sin^2\theta\}\\
&\Leftrightarrow&\ \ (E[ V^{2m_2+1}U_1^{2m_3+1}]-E[ V^{2m_2}U_1^{2m_3+2}] E[VU_1])^2\\
&&<\{E[V^{2m_2}U_1^{2m_3+2}]-(2m_2-1)!!(2m_3-1)!!\}\nonumber\\
&&\ \ \ \ \cdot\{E[V^{2m_2+2}U_1^{2m_3}]+E[ V^{2m_2}U_1^{2m_3+2}](E[VU_1])^2-2E[ V^{2m_2+1}U_1^{2m_3+1}] E[VU_1]\\
&&\ \ \ \ \ \ \ -(2m_2-1)!!(2m_3-1)!!(1-(E[VU_1])^2)\}\\
&\Leftrightarrow&\ \ (E[ V^{2m_2+1}U_1^{2m_3+1}])^2\nonumber\\
&&<E[V^{2m_2}U_1^{2m_3+2}] E[V^{2m_2+2}U_1^{2m_3}]\\
&&\ \ \ -(2m_2-1)!!(2m_3-1)!!\{E[V^{2m_2+2}U_1^{2m_3}]+E[ V^{2m_2}U_1^{2m_3+2}]\\
&&\ \ \ \ \ \ \ \ -2E[ V^{2m_2+1}U_1^{2m_3+1}] E[VU_1]-(2m_2-1)!!(2m_3-1)!!(1-(E[VU_1])^2)\}\\
&\Leftrightarrow&\ \ (E[ V^{2m_2+1}U_1^{2m_3+1}])^2+(2m_2-1)!!(2m_3-1)!!\{E[ V^{2m_2}U_1^{2m_3+2}]+E[V^{2m_2+2}U_1^{2m_3}]\}\\
&&<E[V^{2m_2}U_1^{2m_3+2}] E[V^{2m_2+2}U_1^{2m_3}]+(2m_2-1)!!(2m_3-1)!!\{2E[ V^{2m_2+1}U_1^{2m_3+1}] E[VU_1]\nonumber\\
&&\ \ \ \ \ \ \ +(2m_2-1)!!(2m_3-1)!!(1-(E[VU_1])^2)\}.
\end{eqnarray*}
The proof is complete.\hfill\fbox

The GPI (\ref{jan100}) was established in \cite{LHS} based on the intrinsic connection between moments of Gaussian distributions and the Gaussian hypergeometric functions. Following along these lines, by virtue of Lemma \ref{lem44}, we can show that the GPI (\ref{jan10}) is equivalent to an inequality in terms of  the Gaussian hypergeometric functions (see Lemma \ref{lem55} below).

Denote by $F(a, b; c, z)$ the hypergeometric function (cf. \cite{R}):
$$
    F(a, b; c, z)=\sum_{j=0}^{\infty}\frac{(a)_{(j)}(b)_{(j)}}{(c)_{(j)}}\cdot\frac{z^j}{j!},\quad |z|<1,
$$
where, for $\alpha\not=0$, the factorial function is defined by
$$
    (\alpha)_{(j)}=
   \begin{cases}
   1,\ \  & j=0,\\
\alpha(\alpha+1)\cdots(\alpha+j-1),\ \ & j\ge1.
   \end{cases}
$$

\begin{lem}\label{lem55} Let $m_2,m_3\in \mathbb{N}$. Then, Conjecture \ref{con2} is equivalent to the following claim:

\noindent Claim C:\ \ For $0< |x|<1$,
\begin{eqnarray}\label{qwqw}
&&(2m_2+1)^2\left[F\left(-m_2-1,-m_3;\frac{1}{2},x^2\right)-F\left(-m_2,-m_3;\frac{1}{2},x^2\right)\right]^2\nonumber\\
&<&2
(m_3-m_2)x^2
\left[(2m_2+1)F\left(-m_2-1,-m_3;\frac{1}{2},x^2\right)-1\right]F\left(-m_2,-m_3;\frac{1}{2},x^2\right)\nonumber\\
&&+x^2\left[(2m_2+1)F\left(-m_2,-m_3;\frac{1}{2},x^2\right)-1\right]^2.
\end{eqnarray}
\end{lem}

\noindent{\bf Proof.}\ \  By Lemma \ref{lem44}, we only need to show that Claim B is equivalent to Claim C. Let $(U,V)$ be bivariate Gaussian random variables with $U\sim N(0,1)$ and $V\sim N(0,1)$. Define
$$
x=E[UV].
$$
By the moment formula (cf. \cite{M} and \cite{SL}), we get
\begin{eqnarray}\label{23ZXC}
&&E[ U^{2m_2+1}V^{2m_3+1}]\nonumber\\
&=&\frac{(2m_2+1)!(2m_3+1)!x}{2^{m_2+m_3}}\sum_{j=0}^{\min(m_2,m_3)}\frac{(2x)^{2j}}{(m_2-j)!(m_3-j)!(2j+1)!}\nonumber\\
&=&\frac{(2m_2+1)!(2m_3+1)!x}{m_2!m_3!2^{m_2+m_3}}\sum_{j=0}^{\infty}\frac{(-m_2)_{(j)}(-m_3)_{(j)}}{(\frac{3}{2})_{(j)}}\cdot\frac{(x^2)^{j}}{j!}\nonumber\\
&=&\frac{(2m_2+1)!(2m_3+1)!x}{m_2!m_3!2^{m_2+m_3}}F\left(-m_2,-m_3;\frac{3}{2},x^2\right)\nonumber\\
&=&(2m_2+1)!!(2m_3+1)!!xF\left(-m_2,-m_3;\frac{3}{2},x^2\right),
\end{eqnarray}
\begin{eqnarray*}
&&E[ U^{2m_2}V^{2m_3+2}]\\
&=&\frac{(2m_2)!(2m_3+2)!}{2^{m_2+m_3+1}}\sum_{j=0}^{\min(m_2,m_3+1)}\frac{(2x)^{2j}}{(m_2-j)!(m_3+1-j)!(2j)!}\\
&=&\frac{(2m_2)!(2m_3+2)!}{m_2!(m_3+1)!2^{m_2+m_3+1}}\sum_{j=0}^{\infty}\frac{(-m_2)_{(j)}(-m_3-1)_{(j)}}{(\frac{1}{2})_{(j)}}
\cdot\frac{(x^2)^{j}}{j!}\\
&=&\frac{(2m_2)!(2m_3+2)!}{m_2!(m_3+1)!2^{m_2+m_3+1}}F\left(-m_2,-m_3-1;\frac{1}{2},x^2\right)\nonumber\\
&=&(2m_2-1)!!(2m_3+1)!!F\left(-m_2,-m_3-1;\frac{1}{2},x^2\right),
\end{eqnarray*}
and
\begin{eqnarray}\label{23AA}
&&E[U^{2m_2+2}V^{2m_3}]\nonumber\\
&=&\frac{(2m_2+2)!(2m_3)!}{2^{m_2+m_3+1}}\sum_{j=0}^{\min(m_2+1,m_3)}\frac{(2x)^{2j}}{(m_2+1-j)!(m_3-j)!(2j)!}\nonumber\\
&=&\frac{(2m_2+2)!(2m_3)!}{(m_2+1)!m_3!2^{m_2+m_3+1}}\sum_{j=0}^{\infty}\frac{(-m_2-1)_{(j)}(-m_3)_{(j)}}{(\frac{1}{2})_{(j)}}\cdot\frac{(x^2)^{j}}{j!}
\nonumber\\
&=&\frac{(2m_2+2)!(2m_3)!}{(m_2+1)!m_3!2^{m_2+m_3+1}}F\left(-m_2-1,-m_3;\frac{1}{2},x^2\right)\nonumber\\
&=&(2m_2+1)!!(2m_3-1)!!F\left(-m_2-1,-m_3;\frac{1}{2},x^2\right).
\end{eqnarray}
Then, (\ref{14DD}) is equivalent to
\begin{eqnarray}\label{14E}
&&\left[(2m_2+1)(2m_3+1)xF\left(-m_2,-m_3;\frac{3}{2},x^2\right)\right]^2\nonumber\\
&&+(2m_3+1)F\left(-m_2,-m_3-1;\frac{1}{2},x^2\right)+(2m_2+1)F\left(-m_2-1,-m_3;\frac{1}{2},x^2\right)\nonumber\\
&<&(2m_2+1)(2m_3+1)F\left(-m_2,-m_3-1;\frac{1}{2},x^2\right)F\left(-m_2-1,-m_3;\frac{1}{2},x^2\right)\nonumber\\
&&+2(2m_2+1)(2m_3+1)x^2F\left(-m_2,-m_3;\frac{3}{2},x^2\right)+(1-x^2),\ \ \ \ 0< |x|<1.
\end{eqnarray}

Consider the following functions contiguous to ${F}(a,b;c,z)$:
\begin{equation}\label{contiguous}
  {F}(a\pm 1,b;c,z),\quad {F}(a,b\pm 1; c,z),\quad {F}(a,b;c\pm 1,z).
\end{equation}
To simplify notation,  we denote ${F}(a,b;c,z)$ and
the six  functions in (\ref{contiguous}) respectively by
$$
  {F},\quad {F}(a\pm 1),\quad {F}(b\pm 1),\quad {F}(c\pm 1).
$$
By the following relations of Gauss between contiguous functions (cf. \cite[2.8-(37), (38), page 103]{B}):
\begin{eqnarray*}
&&  c(1-z) F-c F(a-1) +(c-b)z F(c+1)=0,\\
&&
  (b-a)(1-z) F -(c-a)F( a-1) +(c -b) F(b-1)=0,
\end{eqnarray*}
we get
\begin{eqnarray*}
&&x^2F\left(-m_2,-m_3;\frac{3}{2},x^2\right)\\
&=&\frac{1}{2m_3+1}\left[F\left(-m_2-1,-m_3;\frac{1}{2},x^2\right)-(1-x^2)F\left(-m_2,-m_3;\frac{1}{2},x^2\right)\right],
\end{eqnarray*}
and
\begin{eqnarray*}
&&F\left(-m_2,-m_3-1;\frac{1}{2},x^2\right)\\
&=&\frac{1}{2m_3+1}\left[(2m_2+1)F\left(-m_2-1,-m_3;\frac{1}{2},x^2\right)+2(m_3-m_2)(1-x^2)F\left(-m_2,-m_3;\frac{1}{2},x^2\right)\right].
\end{eqnarray*}
Hence, we have that
\begin{eqnarray*}
&&(\ref{14E})\ {\rm holds}\nonumber\\
&\Leftrightarrow&\ \ \frac{(2m_2+1)^2}{x^2}
\left[F\left(-m_2-1,-m_3;\frac{1}{2},x^2\right)-(1-x^2)F\left(-m_2,-m_3;\frac{1}{2},x^2\right)\right]^2\nonumber\\
&&\ \ +2\left[(2m_2+1)F\left(-m_2-1,-m_3;\frac{1}{2},x^2\right)+(m_3-m_2)(1-x^2)F\left(-m_2,-m_3;\frac{1}{2},x^2\right)\right]\nonumber\\
&&<(2m_2+1)
\left[(2m_2+1)F\left(-m_2-1,-m_3;\frac{1}{2},x^2\right)+2(m_3-m_2)(1-x^2)F\left(-m_2,-m_3;\frac{1}{2},x^2\right)\right]\nonumber\\
&&
\ \ \ \cdot F\left(-m_2-1,-m_3;\frac{1}{2},x^2\right)\nonumber\\
&&\ \ +2(2m_2+1)\left[F\left(-m_2-1,-m_3;\frac{1}{2},x^2\right)-(1-x^2)F\left(-m_2,-m_3;\frac{1}{2},x^2\right)\right]+(1-x^2)\nonumber\\
&\Leftrightarrow&\ \ (2m_2+1)^2\left\{\left[F\left(-m_2-1,-m_3;\frac{1}{2},x^2\right)\right]^2
+(1-x^2)\left[F\left(-m_2,-m_3;\frac{1}{2},x^2\right)\right]^2\right\}\nonumber\\
&&\ \ +2(m_2+m_3+1)x^2F\left(-m_2,-m_3;\frac{1}{2},x^2\right)\nonumber\\
&&<2(2m_2+1)
[2m_2+1+(m_3-m_2)x^2]
F\left(-m_2-1,-m_3;\frac{1}{2},x^2\right)F\left(-m_2,-m_3;\frac{1}{2},x^2\right)+x^2\nonumber\\
&\Leftrightarrow&\ \ (2m_2+1)^2\left[F\left(-m_2-1,-m_3;\frac{1}{2},x^2\right)-F\left(-m_2,-m_3;\frac{1}{2},x^2\right)\right]^2\nonumber\\
&&<2
(m_3-m_2)x^2
\left[(2m_2+1)F\left(-m_2-1,-m_3;\frac{1}{2},x^2\right)-1\right]F\left(-m_2,-m_3;\frac{1}{2},x^2\right)\nonumber\\
&&\ \ +x^2\left[(2m_2+1)F\left(-m_2,-m_3;\frac{1}{2},x^2\right)-1\right]^2,\ \ \ \ 0< |x|<1.
\end{eqnarray*}
Therefore, the proof is complete.\hfill\fbox

By Lemmas \ref{lem44} and \ref{lem55}, we find that (\ref{jan10}) $\Leftrightarrow$ (\ref{14DD}) $\Leftrightarrow$ (\ref{qwqw}).  Although we have not been able to prove inequality (\ref{qwqw}) yet, we prove a slightly weaker inequality (see (\ref{tfg}) below). This weaker inequality leads us to derive the following novel moment inequalities for bivariate Gaussian random variables.

\begin{thm}\label{thm3133} Let $m_2,m_3\in\mathbb{N}$ and $(U,V)$ be bivariate Gaussian random variables with $U\sim N(0,1)$ and $V\sim N(0,1)$.

\noindent (i) We have
\begin{eqnarray}\label{Jan55}
E[ U^{2m_2}V^{2m_3}]\le \frac{E[U^{2m_2+2}V^{2m_3}]}{2m_2+1}.
\end{eqnarray}
The equality holds if and only if $U$ and $V$ are independent.

\noindent (ii) Suppose that  $m_3\le m_2$. Then, we have
\begin{eqnarray}\label{23GG}
\left|E[ U^{2m_2+1}V^{2m_3-1}]\right|<\frac{2m_2+1}{2m_3} E[ U^{2m_2}V^{2m_3}].
\end{eqnarray}

\end{thm}

\noindent{\bf Proof.}\ \ First, we prove assertion (i). By \cite[2.8-(20), (21), page 102]{B}, we get
\begin{eqnarray}\label{RFR}
F(a+1)=F+\frac{z}{a} F'=F+\frac{bz}{c} F(a+1,b+1;c+1,z).
\end{eqnarray}
Hence,
\begin{eqnarray}\label{23PPP}
F\left(-m_2-1,-m_3;\frac{1}{2},z\right)=F\left(-m_2,-m_3;\frac{1}{2},z\right)+2m_3zF\left(-m_2,-m_3+1;\frac{3}{2},z\right).
\end{eqnarray}
By (\ref{23AA}) and (\ref{23PPP}), we obtain that
\begin{eqnarray*}
&&E[U^{2m_2+2}V^{2m_3}]-(2m_2+1)E[U^{2m_2}V^{2m_3}]\nonumber\\
&=&\frac{(2m_2+1)!(2m_3)!}{m_2!m_3!2^{m_2+m_3}}F\left(-m_2-1,-m_3;\frac{1}{2},x^2\right)\nonumber\\
&&-\frac{(2m_2+1)!(2m_3)!}{m_2!m_3!2^{m_2+m_3}}F\left(-m_2,-m_3;\frac{1}{2},x^2\right)\nonumber\\
&=&\frac{(2m_2+1)!(2m_3)!}{m_2!(m_3-1)!2^{m_2+m_3-1}}x^2F\left(-m_2,-m_3+1;\frac{3}{2},x^2\right).
\end{eqnarray*}
Then, (\ref{Jan55}) holds.

In the following, we prove assertion (ii). Suppose that
\begin{eqnarray}\label{m2m3}
m_3\le m_2.
\end{eqnarray}

\noindent {\it Step 1.}\ \ For $0\le x\le 1$, define
\begin{eqnarray}\label{explain}
G(x)&:=&(2m_2+1)(1+x)F\left(-m_2,-m_3;\frac{1}{2},x^2\right)\nonumber\\
&&-(2m_2+1)F\left(-m_2-1,-m_3;\frac{1}{2},x^2\right)-x.
\end{eqnarray}
We have
$$
G(0)=0,
$$
and
\begin{eqnarray*}
&&G(1)\\
&=&2(2m_2+1)\frac{\Gamma(\frac{1}{2})\Gamma(m_2+m_3+\frac{1}{2})}{\Gamma(m_2+\frac{1}{2})\Gamma(m_3+\frac{1}{2})}
-(2m_2+1)\frac{\Gamma(\frac{1}{2})\Gamma(m_2+m_3+\frac{3}{2})}{\Gamma(m_2+\frac{3}{2})\Gamma(m_3+\frac{1}{2})}-1\\
&=&(2m_2-2m_3+1)\frac{\Gamma(\frac{1}{2})\Gamma(m_2+m_3+\frac{1}{2})}{\Gamma(m_2+\frac{1}{2})\Gamma(m_3+\frac{1}{2})}-1\\
&\ge&\frac{(m_2+m_3-\frac{1}{2})\cdots(m_3+\frac{1}{2})}{(m_2-\frac{1}{2})\cdots\frac{1}{2}}-1\\
&>&0.
\end{eqnarray*}
We will show that
\begin{eqnarray}\label{23NNN}
G(x)>0,\ \ \ \ 0<x<1.
\end{eqnarray}
Once (\ref{23NNN}) is proved, by (\ref{explain}), we get
\begin{eqnarray}\label{tfg}
&&(2m_2+1)\left[F\left(-m_2-1,-m_3;\frac{1}{2},x^2\right)-F\left(-m_2,-m_3;\frac{1}{2},x^2\right)\right]\nonumber\\
&<&|x|\left[(2m_2+1)F\left(-m_2,-m_3;\frac{1}{2},x^2\right)-1\right],\ \ \ \ 0< |x|<1.
\end{eqnarray}

\noindent {\it Step 2.}\ \ We now prove  (\ref{23NNN}). By (\ref{RFR}), we get
\begin{eqnarray*}
&&G'(x)\\
&=&(2m_2+1)F\left(-m_2,-m_3;\frac{1}{2},x^2\right)\\
&&+\frac{2m_2(2m_2+1)(1+x)}{x}\left[F\left(-m_2,-m_3;\frac{1}{2},x^2\right)-F\left(-m_2+1,-m_3;\frac{1}{2},x^2\right)\right]\\
&&-\frac{2(m_2+1)(2m_2+1)}{x}\left[F\left(-m_2-1,-m_3;\frac{1}{2},x^2\right)-F\left(-m_2,-m_3;\frac{1}{2},x^2\right)\right]-1\\
&=&\frac{(2m_2+1)^2(2+x)}{x}F\left(-m_2,-m_3;\frac{1}{2},x^2\right)\\
&&-\frac{2m_2(2m_2+1)(1+x)}{x}F\left(-m_2+1,-m_3;\frac{1}{2},x^2\right)\\
&&-\frac{2(m_2+1)(2m_2+1)}{x}F\left(-m_2-1,-m_3;\frac{1}{2},x^2\right)-1.
\end{eqnarray*}
Hence
\begin{eqnarray}\label{Jan5A}
&&G'(x)=0\nonumber\\
&\Leftrightarrow&\ \ (2m_2+1)^2(2+x)F\left(-m_2,-m_3;\frac{1}{2},x^2\right)\nonumber\\
&&=2m_2(2m_2+1)(1+x)F\left(-m_2+1,-m_3;\frac{1}{2},x^2\right)\nonumber\\
&&\ \ +2(m_2+1)(2m_2+1)F\left(-m_2-1,-m_3;\frac{1}{2},x^2\right)+x.
\end{eqnarray}
By the following relation of Gauss between contiguous functions (cf. \cite[2.8-(31), page 103]{B}):
$$
[2a-c + (b- a)z]F =a(1 -z)F(a + 1) -(c -a)F(a-1),
$$
we get
\begin{eqnarray*}
&&\left[-2m_2-\frac{1}{2}+(m_2-m_3)x^2\right]F\left(-m_2,-m_3;\frac{1}{2},x^2\right)\nonumber\\
&=&-m_2(1-x^2)F\left(-m_2+1,-m_3;\frac{1}{2},x^2\right)\nonumber\\
&&-\left(m_2+\frac{1}{2}\right)F\left(-m_2-1,-m_3;\frac{1}{2},x^2\right),
\end{eqnarray*}
which implies that
\begin{eqnarray}\label{Jan5B}
&&[(4m_2+1)+2(m_3-m_2)x^2]F\left(-m_2,-m_3;\frac{1}{2},x^2\right)\nonumber\\
&=&2m_2(1-x^2)F\left(-m_2+1,-m_3;\frac{1}{2},x^2\right)\nonumber\\
&&+(2m_2+1)F\left(-m_2-1,-m_3;\frac{1}{2},x^2\right).
\end{eqnarray}
By (\ref{Jan5A}) and (\ref{Jan5B}), we get
\begin{eqnarray*}
&&G'(x)=0\nonumber\\
&\Leftrightarrow&
(2m_2+1)[1-(2m_2+1)x-(2m_3+1)x^2]F\left(-m_2,-m_3;\frac{1}{2},x^2\right)\nonumber\\
&&=(2m_2+1)[1-2(m_2+1)x]F\left(-m_2-1,-m_3;\frac{1}{2},x^2\right)+x(1-x),
\end{eqnarray*}
which implies that
\begin{eqnarray*}
&&G'(x)=0\nonumber\\
&\Leftrightarrow&
F\left(-m_2-1,-m_3;\frac{1}{2},x^2\right)\\
&&=\frac{(2m_2+1)[1-(2m_2+1)x-(2m_3+1)x^2]F\left(-m_2,-m_3;\frac{1}{2},x^2\right)-x(1-x)}{(2m_2+1)[1-2(m_2+1)x]}.
\end{eqnarray*}
Thus, for $x\in(\frac{1}{2(m_2+1)},1)$ with $G'(x)=0$, we have
\begin{eqnarray*}
&&G(x)\\
&=&(2m_2+1)(1+x)F\left(-m_2,-m_3;\frac{1}{2},x^2\right)\\
&&-\frac{(2m_2+1)[1-(2m_2+1)x-(2m_3+1)x^2]F\left(-m_2,-m_3;\frac{1}{2},x^2\right)-x(1-x)}{1-2(m_2+1)x}-x\\
&=&\frac{(2m_2+1)x^2\{[1+2(m_2-m_3)]F\left(-m_2,-m_3;\frac{1}{2},x^2\right)-1\}}{2(m_2+1)x -1}\\
&>&0.
\end{eqnarray*}

By (\ref{23PPP}),  for $0<x<1$,  we have that
\begin{eqnarray}\label{Jan5}
&&\frac{G(x)}{x}\nonumber\\
&=&(2m_2+1)F\left(-m_2,-m_3;\frac{1}{2},x^2\right)-1-2m_3(2m_2+1)xF\left(-m_2,-m_3+1;\frac{3}{2},x^2\right)\nonumber\\
&=&(2m_2+1)m_2!m_3!\sum_{j=0}^{m_2}\frac{(2x)^{2j}}{(m_2-j)!(m_3-j)!(2j)!}-1\nonumber\\
&&-2(2m_2+1)m_2!m_3!x\sum_{j=0}^{m_2-1}\frac{(2x)^{2j}}{(m_2-j)!(m_3-1-j)!(2j+1)!}\nonumber\\
&=&(2m_2+1)-1-2(2m_2+1)m_3x\nonumber\\
&&+(2m_2+1)m_2!m_3!\sum_{j=1}^{m_2}\frac{(2x)^{2j}}{(m_2-j)!(m_3-j)!(2j)!}\nonumber\\
&&-2(2m_2+1)m_2!m_3!x\sum_{j=1}^{m_2-1}\frac{(2x)^{2j}}{(m_2-j)!(m_3-1-j)!(2j+1)!}.
\end{eqnarray}
If $m_2=1$, then $G(x)> 0$ for  $x\in(0, 1)$ by (\ref{Jan5}) and the following elementary inequality:
$$
3|x|<3x^2+1,\ \ \ \ x\in\mathbb{R}.
$$
If $m_2\ge 2$, then $G(x)> 0$ for  $x\in(0, \frac{1}{2(m_2+1)}]$ by (\ref{m2m3}), (\ref{Jan5}) and the following inequality:
\begin{eqnarray*}
&&(2m_2+1)\cdot\frac{1}{m_2+1}-1\\
&&+(2m_2+1)m_2!m_3!\sum_{j=1}^{m_2}\frac{(2x)^{2j}}{(m_2-j)!(m_3-j)!(2j)!}\\
&&-(2m_2+1)m_2!m_3!\sum_{j=1}^{m_2-1}\frac{(2x)^{2j}}{(m_2-j)!(m_3-j)!(2j)!}\nonumber\\
&>&0,\ \ \ \ 0< x<1.
\end{eqnarray*}

Thus far we have shown that $G(0)=0$, $G(x)> 0$ for  $x\in(0, \frac{1}{2(m_2+1)}]$, $G(x)>0$ for $x\in (\frac{1}{2(m_2+1)},1)$ with $G'(x)=0$, and $G(1)>0$.
Then, (\ref{23NNN}) and hence (\ref{tfg}) hold.

\noindent {\it Step 3.}\ \ By (\ref{23PPP}) and (\ref{tfg}), we get
\begin{eqnarray*}
2m_3|x|F\left(-m_2,-m_3+1;\frac{3}{2},x^2\right)<F\left(-m_2,-m_3;\frac{1}{2},x^2\right),\ \ \ \ 0< |x|<1.
\end{eqnarray*}
Thus, by (\ref{23ZXC}) and (\ref{23AA}), we obtain that
$$
\frac{m_2!m_3!2^{m_2+m_3}}{(2m_2+1)!(2m_3-1)!}\left|E[ U^{2m_2+1}V^{2m_3-1}]\right|<\frac{m_2!m_3!2^{m_2+m_3}}{(2m_2)!(2m_3)!} E[ U^{2m_2}V^{2m_3}].
$$
Therefore, (\ref{23GG}) holds.\hfill\fbox

\begin{rem}
To the best of our knowledge, the moment comparison inequality (\ref{23GG}) is not given in the literature. From the proof of Theorem \ref{thm3133}, we can see that (\ref{23GG}) is implied by (\ref{tfg}). By Lemma \ref{lem55}, we find that inequality (\ref{tfg}) is equivalent to the GPI (\ref{jan10}) if $m_3=m_2$ and is weaker than (\ref{jan10}) if $m_3<m_2$.
\end{rem}

\vskip 0.5cm

\begin{large} \noindent\textbf{Acknowledgements} \end{large} We would like to thank the anonymous reviewer for suggestions leading to a significant structural improvement of this paper. This work was supported by the Natural Sciences and Engineering Research Council of Canada (Nos. 559668-2021 and
4394-2018).

\end{document}